\documentclass[12pt]{myarticle}

\newtheorem{prop}{Proposition}
\newtheorem{lemma}{Lemma}
\newtheorem{cor}{Corollary}
\newtheorem{definition}{Definition}
\newtheorem{ex}{Example}
\usepackage{tikz}
\usetikzlibrary{positioning}
\usepackage{subcaption}
\usepackage{caption}
\usepackage{url}
\usepackage{txfonts}
\usepackage[colorlinks]{hyperref}

\newcommand{\n}{\noindent}
\usepackage{graphicx}
\begin{document}

\title{On Neighbourhood Zagreb index of product graphs}
\author[sm]{Sourav Mondal}
\ead{souravmath94@gmail.com}

\address[sm]{Department of mathematics, NIT Durgapur, India.}
\author[nd]{Nilanjan De}
\ead{de.nilanjan@radiffmail.com}

\address[nd]{Department of Basic Sciences and Humanities (Mathematics),\\ Calcutta Institute of Engineering and Management, Kolkata, India.}

\author[ap]{Anita Pal}
\ead{anitabuei@gmail.com}

\address[ap]{Department of mathematics, NIT Durgapur, India.}

\begin{abstract}
There is powerful relation between the chemical behaviour of chemical compounds and their molecular structures. Topological indices defined on these chemical molecular structures are capable to predict physical properties, chemical reactivity and biological activity.
In this article, a new topological index named as Neighbourhood Zagreb index ($M_{N}$) is presented. Here the chemical importance of this newly introduced index is studied and some explicit results for this index of different product graphs such as Cartesian, Tensor and Wreath product is derived. Some of these results are applied to obtain the Neighbourhood Zagreb index of several chemically important graphs and nano-structures.
\medskip

\end{abstract}


\maketitle

\section{Introduction}
 
A chemical graph \cite{tri83,gut86} is a connected graph where loops and parallel edges are not allowed and in which nodes and edges are supposed to be atoms and chemical bonds of compound respectively. Throught this work we use only chemical graphs. Consider a graph $G$  having $V(G)$ and  $E(G)$ as node set and edge set respectively. The degree(valency) of a node $v \in V(G)$, written as $deg_{G}(v)$, is the total number of edges associated with $v$. The set of neighbours of a node $v$ is written as $N_{G}(v)$. For molecular graph, $\vert N_{G}(v) \vert$ = $deg_{G}(v)$.\\
 In mathematical chemistry, molecular descriptors play a leading role specifically in the field of QSPR/QSAR modelling. Amongst them, an outstanding area is preserved for the well-known topological indices or graph invariant. A real valued mapping considering graph as argument is called a graph invariant if it gives same value to graphs which are isomorphic. The order(total count of nodes) and size(total count of edges) of a graph are examples of two graph invariants. In mathematical chemistry, the graph invariants are named as topological indices. Some familiar topological indices are Wiener index, Randi\'c index, connectivity indices, Schultz index, Zagreb indices etc.. The idea of topological indices was initiated when the eminent chemist Harold Wiener found the first topological index, known as Wiener index \cite{das13}, in 1947 for searching boiling points of alkanes. Amidst the topological indices invented on initial stage, the Zagreb indices belong to the well known and well researched molecular descriptors. It was firstly presented by Gutman and Trinajesti\'c \cite{gut72}, where they investigated  how the total energy of $\pi$-electron depends on the structure of molecules and it was recognized on \cite{gut75}. The first ($M_{1}$) and second ($M_{2}$) Zagreb indices are as follows:
 \begin{eqnarray*}
M_{1}(G) = \sum\limits_{u \in V(G)} deg_{G}(u)^{2},
\end{eqnarray*}
 \begin{eqnarray*}
 M_{2}(G) = \sum\limits_{uv \in E(G)}deg_{G}(u) deg_{G}(v).
 \end{eqnarray*}
 $M_{1}(G)$ can also be expressed as 
 \begin{eqnarray*}
  M_{1}(G) = \sum\limits_{uv \in E(G)}[deg_{G}(u) +  deg_{G}(v)].
 \end{eqnarray*}
For more discussion on these indices, we encourage the readers to consult the articles \cite{aza13,ira13,fon14,das13,gut04,ham14,li11}.\\
Let the sum of the degrees of all neighbours of $v$ in $G$ be denoted by $\delta_{G}(v)$, i.e., 
\begin{eqnarray*}
\delta_{G}(v) = \sum_{u \in N_{G}(v)}deg_{G}(u). 
\end{eqnarray*} Following the construction of first Zagreb index, we present here a new degree based topological index named as the Neighbourhood Zagreb index($M_{N}$) which is defined as follows:
 \begin{eqnarray*}
 M_{N}(G) =\sum\limits_{v \in V(G)}\delta_{G}(v)^{2}.
 \end{eqnarray*}
   In mathematical chemistry, graph operations are very significant since certain graphs of chemical interest can be evaluated by various graph operations of different simple graphs. H.Yousefi Azari and co-authors \cite{you08} derived some exact formulae of PI index for Cartesian product of bipartite graphs. P. Paulraja and V.S. Agnes \cite{pau14} evaluated some explicit expressions of the degree distance for the cartesian and wreath products. De et al. \cite{de15} found explicit expressions of the $F$-index under several graph operations. For further illustration on this area, interested readers are suggested some articles \cite{m13,kha09,de14,she17,de16}. We continue this progress for $M_{N}$ index. The objective of this work is to shed some attention on establishing some exact results for the Neighbourhood Zagreb index $M_{N}$ under different product graphs and applying these results for some significant family of chemical graphs and nano-materials.
\section{Chemical significance of the Neighbourhood Zagreb index $(M_{N})$}
\label{sec:1}
According to the instruction of the International Academy of Mathematical Chemistry (IAMC), in order to determine the effectiveness of a  
topological index to predict physiochemical behaviour, we use regression analysis. Usually octane isomers are helpful for such investigation, since they represent a sufficiently large and structurally diverse group of alkanes for the preliminary testing of indices \cite{ran93}. Furtula et al. \cite{fur15} derived that the correlation coefficient of both $M_{1}$ and $F$ for octane isomers is greater than 0.95 with acentric factor and entropy. They also enhanced the skill of prediction of these indices by devising a linear model $(M_{1} + \lambda F)$, where $\lambda$ was varied from -20 to 20. De et al. found that the correlation coefficient of $F$-coindex for octanes in case of  the n-octanes/water partition coefficient (LogP) is 0.966.\\In this article, we find the correlation of entropy (S) and acentric factor (Acent Fac.) with the corresponding Neighbourhood Zagreb index of octanes. The datas of octanes (table:1) are collected from {www.moleculardescriptors.eu/dataset/dataset.htm.} Here we computed that the correlation coefficient ($r$) between acentric factor (Acent Fac.) and $M_{N}$ is -0.99456 (figure.1) and between  entropy (S) and $M_{N}$ is -0.95261 (figure.2). Thus $M_{N}$ can help to predict the Entropy ($r^{2} = 0.98915$) and acentric factor ($r^{2} = 0.90746$) with powerful accuracy. These results confirm the suitability of the indices in QSPR analysis.\\
 A major drawback of most topological indices is their degeneracy, i.e.,two or more isomers possess the same topological index. But this novel index is exceptional for octane isomers. Bonchev et al. \cite{bon81} defined the mean isomer degeneracy as follows :
 \begin{eqnarray*}
  d = \frac{n}{t}
 \end{eqnarray*}
 where $n$ and $t$ are the the number of isomers considered and the number of distinct values that the index assumes for these isomers respectively. Clearly the minimum value of $d$ is 1. As much as '$d$' increases the isomer-discrimination power of topological indices decreases. Thus '$d$' has a significant role for the discriminating power of an index. For octane isomers $M_{N}$ index (table:2) exhibits good response ($d=1$).  
\begin{table}[ht]
\caption{Experimental values of the acentric Factor, entropy(S) and the corresponding values of $M_{N}$.}
\centering
\begin{tabular}{|c| c| c| c| c|}
\hline
\textbf{Octanes} & \textbf{Acent Fac.} & \textbf{S} & \textbf{$M_{N}(G)$} \\ [0.5ex]
\hline
2,2,3,3-tetramethyl butane & 0.255294 & 93.06 & 194 \\
2,3,4-trimethyl pentane & 0.317422 & 102.39 & 144 \\
2,3,3-trimethyl pentane & 0.293177 & 102.06 & 164 \\
2,2,3-trimethyl pentane & 0.300816 & 101.31 & 162 \\
3-methyl-3-ethyl pentane & 0.306899 & 101.48 & 152 \\
2-methyl-3-ethyl pentane & 0.332433 & 106.06 & 132 \\
3,4-dimethyl hexane & 0.340345 & 106.59 & 130 \\
3,3-dimethyl hexane & 0.322596 & 104.74 & 146 \\
2,5-dimethyl hexane & 0.35683 & 105.72 & 118 \\
2,4-dimethyl hexane & 0.344223 & 106.98 & 124 \\
2,3-dimethyl hexane & 0.348247 & 108.02 & 126 \\
2,2-dimethyl hexane & 0.339426 & 103.42 & 138 \\
3-ethyl hexane & 0.362472 & 109.43 & 114 \\
4-methyl heptane & 0.371504 & 109.32 & 110 \\
3-methyl heptane & 0.371002 & 111.26 & 108 \\
2-methyl heptane & 0.377916 & 109.84 & 104 \\
n-octane & 0.397898 & 111.67 & 90 \\
 [1ex]
\hline
\end{tabular}
\label{table:1}
\end{table}
\begin{figure}
\begin{center}
\includegraphics{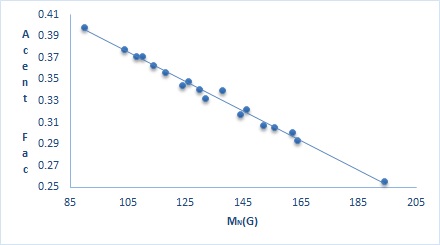}
\caption{Experimental values of Acent Fac. vs. calculated values of $M_{N}$.}
\end{center}
\end{figure}
\begin{figure}
\begin{center}
\includegraphics{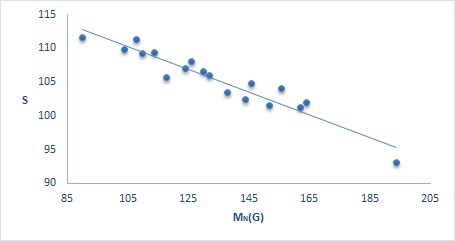}
\caption{Experimental values of entropy(S) vs. calculated values of $M_{N}$.}
\end{center}
\end{figure}
\begin{table}[ht]
\caption{Mean isomer degeneracy ($d$) of different indices for octane isomers.}
\centering
\begin{tabular}{|c| |c|}
\hline
Indices & Mean isomer degeneracy ($d$) \\[1.0ex]
\hline
First Zagreb index ($M_{1}$) & 3.000 \\
\hline
Second Zagreb index ($M_{2}$) & 1.286 \\
\hline
Forgotten topological index ($F$) & 2.571 \\
\hline
Hosoya index ($Z$) & 1.286 \\
\hline
Merrifield-Simmons index($\sigma$) & 1.200 \\
\hline
Connectivity index($\chi$) & 1.125 \\
\hline
Harary index($\eta$) & 1.059 \\
\hline
\textbf{Neighbourhood Zagreb index ($M_{N}$)} & \textbf{1.000}\\
[1ex]
\hline
\end{tabular}
\label{table:2}
\end{table}
\section{Main Result}
In this section, we evaluate the newly introduced index of different product graphs such as Cartesian, Wreath and Tensor product of graphs. Here $V_{i}$ and $E_{i}$ denote vertex set and edge set of $G_{i}$ respectively. Also for path, cycle and complete graphs with $n$ nodes, we use $P_{n}$, $C_{n}$ and $K_{n}$ respectively. We proceed with the following lemma directly followed from definitions.
\begin{lemma}{\label{lem1}} For graph $G$,  we have
\begin{enumerate}
\item[(i)]
$\sum\limits_{v \in V(G)} \delta_{G} (v) = M_{1}(G)$,
\item[(ii)]
$\sum\limits_{v \in V(G)} deg_{G} (v) \delta_{G} (v) = 2M_{2}(G)$.
\end{enumerate}
\end{lemma} 
\subsection{Cartesian product}
\begin{definition}
The Cartesian product of $G_{1}$, $G_{2}$, written as $ G_{1}$ $\times$ $G_{2}$, containing node set $V_{1}$ $\times$ $V_{2} $ and $(u_{1},v_{1})$ is adjacent with $ (u_{2},v_{2}) $  iff [$ u_{1}=u_{2} $ and $ v_{1}v_{2}$ $\in$ $E_{2} $] or [$ v_{1}=v_{2} $ and $ u_{1}u_{2}$ $\in$ $E_{1} $].   
\end{definition}
Clearly the above definition yield the lemma stated below.
\begin{lemma}{\label{lem2}} For graphs $G_{1}$ and $G_{2}$, we have
\begin{enumerate}
\item[(i)]
  $\delta _{G_{1} \times G_{2}} (u,v) = \delta_{G_{1}}(u)+\delta_{G_{2}}(v)+2deg_{G_{1}}(u)deg_{G_{2}}(v)$,
 \item[(ii)]
 $\vert E(G_{1} \times G_{2})\vert = \vert V_{2} \vert \vert E_{1}\vert  + \vert V_{1} \vert \vert E_{2} \vert$.
\end{enumerate}
\end{lemma}
 In \cite{kla96,kha08} different topological indices for Cartesian product of graphs were studied. Here we move to calculate the Neighbourhood Zagreb index of Cartesian product graphs.
\begin{prop}{\label{pro1}} The Neighbourhood Zagreb index of $G_{1} \times G_{2}$ is given by
\begin{eqnarray*}
M_{N}(G_{1} \times G_{2}) &=& 6M_{1}(G_{1})M_{1}G_{2}) + \vert V_{2} \vert M_{N}(G_{1}) + \vert V_{1} \vert M_{N}(G_{2}) + 16[\vert E_{2} \vert M_{2}(G_{1})\\
&&  + \vert E_{1} \vert M_{2}(G_{2})].
\end{eqnarray*}
\end{prop}
\n\textit{Proof}. From definition of Neighbourhood Zagreb index and applying lemma \ref{lem2} and lemma \ref{lem1}, we get
\begin{eqnarray*}
M_{N}(G_{1} \times G_{2})&=&\sum\limits_{(u,v) \in V_{1} \times V_{2}}\delta_{G_{1} \times G_{2}}^{2}(u,v)\\\\&=&\sum\limits_{u \in V_{1}}\sum\limits_{v \in V_{2}}[\delta_{G_{1}}(u) + \delta_{G_{2}}(v) + 2deg_{G_{1}}(u)deg_{G_{2}}(v)]^{2}\\\\&=&6M_{1}(G_{1})M_{1}G_{2}) + \vert V_{2} \vert M_{N}(G_{1}) + \vert V_{1} \vert M_{N}(G_{2}) + 16[\vert E_{2} \vert M_{2}(G_{1})\\ 
&&+ \vert E_{1} \vert M_{2}(G_{2})].
\end{eqnarray*}
Hence the result.                  
Using the proposition \ref{pro1}, we have the following results.
\begin{ex} 
The Cartesian product of $P_{2}$ and $P_{n+1}$ produces the ladder graph $L_{n}$ (figure.3). By the above proposition, we derive the following result.
\begin{eqnarray*}
M_{N}(L_{n}) &=& 162n-132.
\end{eqnarray*}
\end{ex}
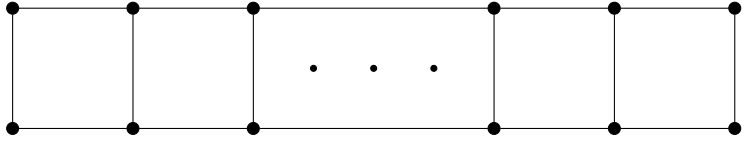
\begin{figure*}[h]
\begin{center}
\begin{tikzpicture}[scale=.4]
\tikzstyle{every node}=[draw, shape=circle, fill=black, scale=.4]
  \node (n1) at (1,1) {};
  \node (n2) at (5,1)  {};
  \node (n3) at (9,1)  {};
  \node (n4) at (17,1) {};
  \node (n5) at (21,1)  {};
  \node (n6) at (25,1)  {};
  \node (n7) at (1,5)  {};
  \node (n8) at (5,5)  {};
  \node (n9) at (9,5)   {};
  \node (n10) at (17,5)   {};
  \node (n11) at (21,5)   {};
  \node (n12) at (25,5)   {};
\tikzstyle{every node}=[draw, shape=circle, fill=black, scale=.2]
  \node  (n13) at (11,3)  {};
  \node   (n14) at (13,3) {};
  \node  (n15)  at (15,3) {};
  \foreach \from/\to in {n1/n2,n2/n3,n3/n4,n4/n5,n5/n6,n7/n8,n8/n9,n9/n10,n10/n11,n11/n12,n1/n7,n2/n8,n3/n9,n4/n10,n5/n11,n6/n12}
    \draw (\from) -- (\to);
\end{tikzpicture}
\caption{The ladder graph $L_n$.}
\end{center}
\end{figure*}
\begin{ex}
For a $C_{4}-nanotorus$  $TC_{4}(m,n) = C_{m} \times C_{n}$ , the Neighbourhood Zagreb index is given by 
\begin{eqnarray*}
M_{N}(TC_{4}(m,n)) &=& 256mn.
\end{eqnarray*}
\end{ex}
\begin{ex}
The Cartesian product of $P_{n}$ and $C_{4}$ yields a $C_{4}-nanotube$ $TUC_{4}(m,n) = P_{n} \times C_{m}$. Its $M_{N}$ index is as follows:
\begin{eqnarray*}
M_{N}(TUC_{4}(m,n)) &=& 256mn-374m,   n\geq4. 
\end{eqnarray*}
\end{ex}
 \begin{ex}
 The Neighbourhood Zagreb index of the grids $(P_{n} \times P_{m})$ is given by
 \begin{eqnarray*}
 M_{N}(P_{n} \times P_{m}) &=& 256mn - 310m - 310n + 216, m,n\geq4.
\end{eqnarray*}
\end{ex}
\begin{ex}
For a $n$-prism($K_{2} \times C_{n}$) (figure.4), the Neighbourhood Zagreb index is given below. 
\begin{eqnarray*}
M_{N}(K_{2} \times C_{n}) &=& 162n.
\end{eqnarray*}
\end{ex}
\begin{ex} 
The Cartesian product of $K_{n}$ and $K_{m}$ yields the Rook's graph (figure.4). Applying the above proposition we have the following.
\begin{eqnarray*}
M_{N}(K_{m} \times K_{n}) &=& mn[6(m-1)^{2}(n-1)^{2} + (n-1)^{4} + (m-1)^{4}\\ 
&&+ 4(m-1)(n-1)\lbrace (m-1)^{2} + (n-1)^{2} \rbrace].
\end{eqnarray*}
\end{ex}
\begin{figure*}[h]
\begin{center}
\captionsetup{singlelinecheck = true}
\begin{tikzpicture}[scale=.4]
\tikzstyle{every node}=[draw, shape=circle, fill=black, scale=.4]
  \node (n1) at (6,6) {};
  \node (n2) at (3,1)  {};
  \node (n3) at (9,1)  {};
  \node (n4) at (15,8) {};
  \node (n5) at (12,3)  {};
  \node (n6) at (18,3)  {};  
  \foreach \from/\to in {n1/n2,n2/n3,n3/n1,n4/n5,n5/n6,n6/n4,n1/n5,n6/n1,n4/n2,n6/n2,n3/n4,n3/n5,n1/n4,n2/n5,n3/n6}
    \draw (\from) -- (\to);
\end{tikzpicture}
\hspace{2cm}
\begin{tikzpicture}[scale=.4]
\tikzstyle{every node}=[draw, shape=circle, fill=black, scale=.4]
  \node (n1) at (1,4) {};
  \node (n2) at (2.5,1)  {};
  \node (n3) at (6.5,1)  {};
  \node (n4) at (8,4) {};
  \node (n5) at (6.5,7)  {};
  \node (n6) at (2.5,7)  {};
  \node (n7) at (2.5,4)  {};
  \node (n8) at (3.5,2)  {};
  \node (n9) at (5.5,2)   {};
  \node (n10) at (6.5,4)   {};
  \node (n11) at (5.5,6)   {};
  \node (n12) at (3.5,6)   {};
  \foreach \from/\to in {n1/n2,n2/n3,n3/n4,n4/n5,n5/n6,n6/n1,n7/n8,n8/n9,n9/n10,n10/n11,n11/n12,n12/n7,n1/n7,n2/n8,n3/n9,n4/n10,n5/n11,n6/n12}
    \draw (\from) -- (\to);
\end{tikzpicture}
\caption*{$K_{3} \times K_{2}$ \hspace{7cm} $K_{2} \times C_{6}$}
\captionsetup{singlelinecheck = true, justification=justified}
\caption{The example of Rook's graph ($K_{3} \times K_{2}$) and $n$-Prism graph ($n=6$).}
\end{center}
\end{figure*}
 Now we generalize the Proposition \ref{pro1}. We begin with the following lemma.
\begin{lemma}{\label{lem3}}If $G_{1}$, $G_{2}$,....., $G_{n}$ be $n$ graphs and $ V = V(\bigotimes\limits_{i=1}^{n} G_{i})$, $ E = E(\bigotimes\limits_{i=1}^{n} G_{i})$, then we have
\begin{enumerate}
\item[(i)]
$\vert E(\bigotimes\limits_{i=1}^{n}G_{i})\vert$ = $\vert V \vert \sum\limits_{i=1}^{n}\frac{\vert E_{i} \vert}{\vert V_{i} \vert}$,
\item[(ii)]
$M_{1}(\bigotimes\limits_{i=1}^{n}G_{i}) = \vert V \vert \sum\limits_{i=1}^{n} \frac{M_{1}(G_{i})}{\vert V_{i} \vert} + 4\vert V \vert \sum\limits_{i \neq j, i,j = 1}^{n-1} \frac{\vert E_{i} \vert \vert E_{j} \vert}{\vert V_{i} \vert \vert V_{j} \vert}$,
\item[(iii)]
$M_{2}(\bigotimes\limits_{i=1}^{n}G_{i}) = \vert V \vert \sum\limits_{i=1}^{n} \frac{M_{2}(G_{i})}{\vert V_{i} \vert} + 3 \sum\limits_{i=1}^{n} M_{1}(G_{i}) ( \frac{\vert E \vert}{\vert V_{i} \vert}-\frac{\vert V\vert \vert E_{i} \vert}{\vert V_{i} \vert^{2}}) + 4\vert V \vert \sum\limits_{i \neq j \neq k , i,j,k = 1}^{n} \frac{\vert E_{i} \vert \vert E_{j} \vert \vert   E_{k} \vert}{\vert V_{i} \vert \vert V_{j} \vert \vert V_{k} \vert}$.
\end{enumerate}
\end{lemma}
 \n\textit{Proof}. Applying lemma \ref{lem2}(ii) and an inductive argument, (i) is clear. In order to proof (ii) and (iii), we refer to \cite{kha09}.           
\begin{prop}{\label{pro2}} If  $G_{1}$, $G_{2}$,..., $G_{n}$ be $n$ graphs, then we have
 \begin{eqnarray*}
M_{N}(\bigotimes\limits_{i=1}^{n}G_{i}) &=&  \vert V \vert  \sum\limits_{i=1}^{n} \frac{M_{N}(G_{i})}{\vert V_{i} \vert} + 3 \vert V \vert \sum\limits_{i \neq j, i,j = 1 }^{n}\frac{M_{1}(G_{i}) M_{1}(G_{j})}{\vert V_{i} \vert   \vert V_{j} \vert}\\
&& + 24 \vert V \vert \sum\limits_{i \neq j \neq k, i,j,k = 1}^{n}\frac{M_{1}(G_{i})\vert E_{j} \vert \vert E_{k} \vert}{\vert V_{i} \vert \vert V_{j} \vert \vert V_{k} \vert} + 16 \vert V \vert \sum\limits_{i \neq j, i,j = 1}^{n} \frac{M_{2}(G_{i}) \vert E_{j} \vert }{\vert V_{i} \vert   \vert V_{j} \vert} \\
&&+ 16 \vert V \vert \sum\limits_{i \neq j \neq k \neq l, i,j,k,l = 1}^{n} \frac{\vert E_{i} \vert \vert E_{j} \vert \vert   E_{k} \vert \vert E_{l} \vert}{\vert V_{i} \vert \vert V_{j} \vert \vert V_{k} \vert \vert V_{l} \vert}.        
\end{eqnarray*}
\end{prop}
\n\textit{Proof}. We derive the formula by mathematical induction. Evidently the result holds for $n = 2$. Let us take the proposition to be true for ($n-1$)  graphs. Then we obtain 
\begin{eqnarray*}
 M_{N}(\bigotimes\limits_{i=1}^{n}G_{i}) &=&  M_{N}(\bigotimes_{i = 1}^{n-1}G_{i} \times G_{n})\\
 &&= 6 M_{1}(\bigotimes\limits_{i=1}^{n-1}G_{i}) M_{1}(G_{n}) + \vert V_{n} \vert M_{N}(\bigotimes\limits_{i=1}^{n-1}G_{i}) + \vert V(\bigotimes\limits_{i=1}^{n-1}G_{i}) \vert  M_{N}(G_{n})\\
 &&+16[ M_{2}(\bigotimes\limits_{i=1}^{n-1}G_{i}) \vert E_{n} \vert + M_{2}(G_{n}) \vert E(\bigotimes\limits_{i=1}^{n-1}G_{i}) \vert ].
 \end{eqnarray*}
 Now applying lemma \ref{lem3}, we get
 \begin{eqnarray*}  
M_{N}(\bigotimes\limits_{i=1}^{n}G_{i}) &=& 6 \vert V \vert \frac{M_{1}(G_{n})}{\vert V_{n} \vert}[ \sum\limits_{i=1}^{n-1} \frac{M_{1}(G_{i})}{\vert V_{i} \vert} + 4 \sum\limits_{i \neq j, i,j = 1}^{n-1} \frac{\vert E_{i} \vert \vert E_{j} \vert}{\vert V_{i} \vert \vert V_{j} \vert}
\end{eqnarray*}
\begin{eqnarray*}
 &&+ \vert V \vert [  \sum\limits_{i=1}^{n-1} \frac{M_{N}(G_{i})}{\vert V_{i} \vert} + 3 \sum\limits_{i \neq j, i,j = 1 }^{n-1}\frac{M_{1}(G_{i}) M_{1}(G_{j})}{\vert V_{i} \vert   \vert V_{j} \vert}\\
     &&+ 24 \sum\limits_{i \neq j \neq k, i,j,k = 1}^{n-1}\frac{M_{1}(G_{i}) \vert E_{j} \vert \vert E_{k} \vert}{\vert V_{i} \vert \vert V_{j} \vert \vert V_{k} \vert} + 16 \sum\limits_{i \neq j, i,j = 1}^{n-1} \frac{M_{2}(G_{i}) \vert E_{j} \vert }{\vert V_{i} \vert   \vert V_{j} \vert}\\
     &&+ 16 \sum\limits_{i \neq j \neq k \neq l, i,j,k,l = 1}^{n} \frac{\vert E_{i} \vert \vert E_{j} \vert \vert   E_{k} \vert \vert E_{l} \vert}{\vert V_{i} \vert \vert V_{j} \vert \vert V_{k} \vert \vert V_{l} \vert} ] + \vert V \vert \frac{M_{N}(G_{n})}{\vert V_{n} \vert}\\
     && + 16 \frac{\vert V \vert}{\vert V_{n} \vert} [\vert E_{n} \vert \lbrace \sum\limits_{i=1}^{n-1} \frac{M_{2}(G_{i})}{\vert V_{i} \vert} + 3 \sum\limits_{i=1}^{n-1} \frac{M_{1}(G_{i})}{\vert V_{i} \vert} ( \sum\limits_{j=1}^{n-1} \frac{\vert E_{j} \vert}{\vert V_{j} \vert}-\frac{\vert E_{i} \vert}{\vert V_{i} \vert})\\
     &&+ 4 \sum\limits_{i \neq j \neq k , i,j,k = 1}^{n-1} \frac{\vert E_{i} \vert \vert E_{j} \vert \vert   E_{k} \vert}{\vert V_{i} \vert \vert V_{j} \vert \vert V_{k} \vert} \rbrace + M_{2}(G_{n}) \sum\limits_{i=1}^{n-1} \frac{\vert E_{i} \vert}{\vert V_{i} \vert}.
  \end{eqnarray*}
After simplification, we have
 \begin{eqnarray*}
  M_{N}(\bigotimes\limits_{i=1}^{n} G_{i}) &=& \vert V \vert \sum\limits_{i=1}^{n} \frac{M_{N}(G_{i})}{\vert V_{i} \vert} + 3 \vert V \vert [\sum\limits_{i \neq j, i,j = 1 }^{n-1}\frac{M_{1}(G_{i}) M_{1}(G_{j})}{\vert V_{i} \vert   \vert V_{j} \vert} \\
  &&+ 2 \frac{M_{1}(G_{n})}{\vert V_{n} \vert} \sum\limits_{i=1}^{n-1} \frac{M_{N}(G_{i})}{\vert V_{i} \vert}]+ 24 \vert V \vert [ \sum\limits_{i \neq j \neq k, i,j,k = 1}^{n-1}\frac{M_{1}(G_{i}) \vert E_{j} \vert \vert E_{k} \vert}{\vert V_{i} \vert \vert V_{j} \vert \vert V_{k} \vert}\\
  && + \sum\limits_{i \neq j, i,j = 1}^{n-1}\frac{M_{1}(G_{n}) \vert E_{i} \vert \vert E_{j} \vert}{\vert V_{n} \vert \vert V_{i} \vert \vert V_{j} \vert} + 2\lbrace \sum\limits_{i,j = 1}^{n-1}\frac{M_{1}(G_{i}) \vert E_{j} \vert \vert E_{n} \vert}{\vert V_{i} \vert \vert V_{j} \vert \vert V_{n} \vert} \\
  &&- \sum\limits_{i = 1}^{n-1}\frac{M_{1}(G_{i}) \vert E_{i} \vert \vert E_{n} \vert}{\vert V_{i} \vert \vert V_{i} \vert \vert V_{n} \vert}\rbrace] + 16 \vert V \vert [ \sum\limits_{i \neq j, i,j = 1}^{n-1} \frac{M_{2}(G_{i}) \vert E_{j} \vert }{\vert V_{i} \vert   \vert V_{j} \vert} \\
  &&+ \sum\limits_{ i = 1}^{n-1} \frac{M_{2}(G_{n}) \vert E_{i} \vert }{\vert V_{n} \vert   \vert V_{i} \vert} + \sum\limits_{i = 1}^{n-1} \frac{M_{2}(G_{i}) \vert E_{n} \vert }{\vert V_{i} \vert   \vert V_{n} \vert}]\\
  &&+ 16 \vert V \vert [ \sum\limits_{i \neq j \neq k \neq l, i,j,k,l = 1}^{n-1} \frac{\vert E_{i} \vert \vert E_{j} \vert \vert   E_{k} \vert \vert E_{l} \vert}{\vert V_{i} \vert \vert V_{j} \vert \vert V_{k} \vert \vert V_{l} \vert}\\
  &&+ 4 \sum\limits_{i \neq j \neq k \neq l, i,j,k,l = 1}^{n-1} \frac{\vert E_{i} \vert \vert E_{j} \vert \vert   E_{k} \vert \vert E_{n} \vert}{\vert V_{i} \vert \vert V_{j} \vert \vert V_{k} \vert \vert V_{n} \vert}]\\
     && =  \vert V \vert  \sum\limits_{i=1}^{n} \frac{M_{N}(G_{i})}{\vert V_{i} \vert} + 3 \vert V \vert \sum\limits_{i \neq j, i,j = 1 }^{n}\frac{M_{1}(G_{i}) M_{1}(G_{j})}{\vert V_{i} \vert   \vert V_{j} \vert}\\
     && + 24 \vert V \vert \sum\limits_{i \neq j \neq k, i,j,k = 1}^{n}\frac{M_{1}(G_{i}) \vert E_{j} \vert \vert E_{k} \vert}{\vert V_{i} \vert \vert V_{j} \vert \vert V_{k} \vert}+ 16 \vert V \vert \sum\limits_{i \neq j, i,j = 1}^{n} \frac{M_{2}(G_{i}) \vert E_{j} \vert }{\vert V_{i} \vert   \vert V_{j} \vert}\\
     && + 16 \vert V \vert \sum\limits_{i \neq j \neq k \neq l, i,j,k,l = 1}^{n} \frac{\vert E_{i} \vert \vert E_{j} \vert \vert   E_{k} \vert \vert E_{l} \vert}{\vert V_{i} \vert \vert V_{j} \vert \vert V_{k} \vert \vert V_{l} \vert}.                                                               
\end{eqnarray*}
Hence the derived result follows.               
\begin{definition} Consider the graph $G$ containing $m$-tuples $b_{1}, b_{2}, ..., b_{m}$ with $b_{i}$ $\in \lbrace 0, 1, ..., n_{i}-1 \rbrace, n_{i}\geqslant 2$, as vertices and let whenever the difference of two tuples is exactly one place, the corresponding two vertices are adjacent. Such a graph is known as Hamming graph. 
The necessery and sufficient criteria for a graph $G$ to be a Hamming graph is that $G = \bigotimes\limits_{i = 1}^{m} K_{n_{i}}$ and that is why such a graph $G$ is naturally written as $H_{n_{1},n_{2},...,n_{m}}$.
\end{definition}
Hamming graph is very useful in coding theory specially in error correcting codes. Also such type of graph is effective in association schemes. Applying proposition \ref{pro2}, we have the corollary stated below.
\begin{cor} The Neighbourhood Zagreb index of Hamming graph is given by 
\begin{eqnarray*}
M_{N}(G) &=& \prod\limits_{i=1}^{m}n_{i}[\sum\limits_{i=1}^{m} (n_{i}-1)^{4} + 3 \sum\limits_{i \neq j, i,j=1}^{m} (n_{i}-1)^{2}(n_{j}-1)^{2} \\
&&+ 6\sum\limits_{i \neq j \neq k, i,j,k=1}^{m} (n_{i}-1)^{2}(n_{j}-1)(n_{k}-1) +4\sum\limits_{i \neq j, i,j=1}^{m} (n_{i}-1)^{3}(n_{j}-1) \\
&&+ \sum\limits_{i \neq j \neq k \neq l, i,j,k,l=1}^{m} (n_{i}-1)(n_{j}-1)(n_{k}-1)(n_{l}-1)].
\end{eqnarray*}
\end{cor}
\begin{ex}
When $n_{1}$, $n_{2}$, ....,$ n_{m}$ are all equal to 2, the hamming graph $G$ is  known as a hypercube (figure.5) with dimension $m$ and written as $Q_{m}$. We compute the following.
\begin{eqnarray*}
M_{N}(Q_{m}) &=& 2^{m}m^{4}.
\end{eqnarray*}
\begin{figure*}[h]
\begin{center}
\captionsetup{singlelinecheck = true}
\begin{tikzpicture}[scale=.4]
\tikzstyle{every node}=[draw, shape=circle, fill=black, scale=.4]
  \node (n1) at (1,1) {};
  \node (n2) at (5,1)  {};
  \node (n3) at (7,3)  {};
  \node (n4) at (3,3) {};
  \node (n5) at (1,5)  {};
  \node (n6) at (5,5)  {};
  \node (n7) at (7,7)  {};
  \node (n8) at (3,7)  {};
  \node (n9) at (9,1) {};
  \node (n10) at (13,1)  {};
  \node (n11) at (15,3)  {};
  \node (n12) at (11,3) {};
  \node (n13) at (9,5)  {};
  \node (n14) at (13,5)  {};
  \node (n15) at (15,7)  {};
  \node (n16) at (11,7)  {};
  \foreach \from/\to in {n1/n2,n2/n3,n3/n4,n4/n1,n5/n6,n6/n7,n7/n8,n8/n5,n1/n5,n2/n6,n3/n7,n4/n8,n9/n10,n10/n11,n11/n12,n12/n9,n13/n14,n14/n15,n15/n16,n16/n13,n9/n13,n10/n14,n11/n15,n12/n16}
    \draw (\from) -- (\to);
    \draw[-] (n1) to [bend right=30] (n9);
    \draw[-] (n2) to [bend right=30] (n10);
    \draw[-] (n3) to [bend right=30] (n11);
    \draw[-] (n4) to [bend right=30] (n12);
    \draw[-] (n5) to [bend left=30] (n13);
    \draw[-] (n6) to [bend left=30] (n14);
    \draw[-] (n7) to [bend left=30] (n15);
    \draw[-] (n8) to [bend left=30] (n16);
\end{tikzpicture}
\caption*{$Q_{4}$}
\captionsetup{singlelinecheck = true, justification=justified}
\caption{Example of Hypercube.}
\end{center}
\end{figure*}
\end{ex}
\subsection{Tensor product}
\begin{definition}
The tensor product of $G_{1}$, $G_{2}$, written as $G_{1} \circledast G_{2}$, is a graph containing node set $V_{1} \times V_{2}$ and $(u_{1},v_{1})$ is adjacent with $(u_{2}, v_{2})$ iff $u_{1}u_{2} \in E_{1}$ and $v_{1}v_{2} \in E_{2}$.
\end{definition}
Clearly the definition gives the lemma as follows:
\begin{lemma}{\label{lem4}} For graphs $G_{1}$ and $G_{2}$, we have 
\begin{eqnarray*}
\delta_{G_{1} \circledast G_{2}}(u,v)&=&\delta_{G_{1}}(u)\delta_{G_{2}}(v).
\end{eqnarray*}
\end{lemma}
Z. Yarahmadi discussed about the Randi\'c, arithmetic,  geometric and Zagreb indices of tensor product of two graphs in \cite{yar11}. Also in \cite{pat12,na1} various topological descriptors of Tensor product graphs are calculated. Here we proceed for the $M_{N}$ index of tensor product of two graphs. 
\begin{prop}{\label{pro3}} The Neighbourhood Zagreb index of tensor product for two graphs is given by
\begin{eqnarray*}
M_{N}(G_{1} \circledast G_{2})&=& M_{N}(G_{1})M_{N}(G_{2}).
\end{eqnarray*}
\end{prop}
\n\textit{Proof}. By the definition of the $M_{N}$ index and applying lemma \ref{lem4}, we get 
\begin{eqnarray*}
M_{N}(G_{1} \circledast G_{2})&=& \sum\limits_{(v_{1},v_{2}) \in V_{1} \times V_{2}}\delta_{G_{1} \times G_{2}}^{2}(v_{1},v_{2})\\\\&=&\sum\limits_{v_{1} \in V_{1}}\sum\limits_{v_{2} \in V_{2}}[\delta_{G_{1}}(v_{1})\delta_{G_{2}}(v_{2})]^{2}\\\\&=&M_{N}(G_{1})M_{N}(G_{2}).
\end{eqnarray*}
Which is the required result.        
\begin{ex}
Usinng the proposition \ref{pro3}, we have the following computations.
\begin{enumerate}
\item[(i)]
$M_{N}$ ($P_{n} \circledast P_{m}$) = $(16n-38)(16m-38)$, $m,n \geqslant 4$.
\item[(ii)]
$M_{N}$($C_{n} \circledast C_{m}$) = $256mn$.
\item[(iii)]
$M_{N}$($K_{n} \circledast K_{m}$)= $mn(m-1)^{4}(n-1)^{4}$.
\item[(iv)]
$M_{N}$($P_{n} \circledast C_{m}$)= $16m(16n-38$), $n \geqslant 4$.
\item[(v)]
$M_{N}$($P_{n} \circledast K_{m}$)= $m(m-1)^{4}(16n-38)$, $n \geqslant 4$.
\item[(vi)]
$M_{N}$($C_{n} \circledast K_{m}$) = $16mn(m-1)^{4}$.
\end{enumerate}
\end{ex}
\subsection{Wreath Product}
\begin{definition}
 The Wreath product (also known as composition) of $G_{1}$ and $G_{2}$ having  $V_{1}$ and $V_{2}$ as vertex sets with no common vertex and edge sets $ E_{1} $ and $ E_{2} $ is the graph $ G_{1} [G_{2}]$  containing node set $V_{1} \times V_{2}$ and $(u_{1},v_{1})$ is adjacent to $(u_{2},v_{2})$ iff ($u_{1}u_{2} \in E_{1})$ or ($u_{1}=u_{2}$ and $v_{1}v_{2} \in E_{2}$).
 \end{definition}
  The following lemma is obvious from the definition. 
 \begin{lemma}{\label{lem5}} For graphs $G_{1}$ and $G_{2}$ , we have
  \begin{eqnarray*}
 \delta_{G_{1}[G_{2}]}(u,v) &=& \vert V_{2} \vert^{2}\delta_{G_{1}}(u) + \delta_{G_{2}}(v) + 2\vert E_{2} \vert deg_{G_{1}}(u) + \vert V_{2} \vert deg_{G_{1}}(u)deg_{G_{2}}(v).
 \end{eqnarray*}
 \end{lemma}
 In \cite{de15,de16,don17} different topological indices for Wreath product of graphs are derived. Here we preserve this movement to find the $M_{N}$ index of Wreath product of two graphs.
 \begin{prop}{\label{pro4}} The Neighbourhood Zagreb index of Wreath product for two graphs is given by
 \begin{eqnarray*}
 M_{N}(G_{1}[G_{2}])&=&\vert V_{2} \vert^{4}M_{N}(G_{1}) + M_{N}(G_{2}) + 12\vert V_{2} \vert \vert E_{2} \vert^{2} M_{1}(G_{1}) + 8\vert E_{1} \vert \vert E_{2} \vert M_{1}(G_{2})\\
 && + 8\vert V_{2} \vert^{2} \vert E_{2} \vert (\vert V_{2} \vert + 1)M_{2}(G_{1}) + 8\vert V_{2} \vert \vert E_{1} \vert M_{2}(G_{2})\\
 && + 3\vert V_{2} \vert^{2}M_{1}(G_{1})M_{1}(G_{2}).
 \end{eqnarray*}
 \end{prop}
 \n\textit{Proof}. From definition of Neighbourhood Zagreb index and using lemma \ref{lem5}, we have 
 \begin{eqnarray*}
 M_{N}(G_{1}[G_{2}])&=&\sum\limits_{(u,v) \in V_{1} \times V_{2}} \delta_{G_{1}[G_{2}]}^{2}(u,v)\\&=&\sum\limits_{u \in V_{1}}\sum\limits_{v \in V_{2}}[\vert V_{2} \vert^{2}\delta_{G_{1}}(u) + \delta_{G_{2}}(v) + 2\vert E_{2} \vert deg_{G_{1}}(u)\\
 && + \vert V_{2} \vert deg_{G_{1}}(u)deg_{G_{2}}(v)]^{2}\\&=& \sum\limits_{u \in V_{1}}\sum\limits_{v \in V_{2}}[\vert V_{2} \vert^{4} \delta_{G_{1}}(u)^{2} + \vert V_{2} \vert^{2} deg_{G_{1}}(u)^{2} deg_{G_{2}}(v)^{2} + 4\vert E_{2} \vert^{2} deg_{G_{1}}(u)^{2}\\
 && + \delta_{G_{2}}(v)^{2}+ 2 \vert V_{2} \vert^{3} \delta_{G_{1}}(u)deg_{G_{1}}(u)deg_{G_{2}}(v)\\
 && + 4 \vert V_{2} \vert^{2} \vert E_{2} \vert \delta_{G_{1}}(u)deg_{G_{1}}(u)+ 2\vert V_{2} \vert^{2} \delta_{G_{1}}(u) \delta_{G_{2}}(v)\\
 &&+ 4 \vert V_{2} \vert \vert E_{2} \vert deg_{G_{1}}(u)^{2} deg_{G_{2}}(v) + 2 \vert V_{2} \vert \delta_{G_{2}}(v) deg_{G_{1}}(u) deg_{G_{2}}(v)\\
 && + 4 deg_{G_{1}}(u) \vert E_{2} \vert \delta_{G_{2}}(v)]. 
 \end{eqnarray*}
 Applying lemma \ref{lem1}, we have
 \begin{eqnarray*}
M_{N}(G_{1}[G_{2}])&=&\vert V_{2} \vert^{4}M_{N}(G_{1}) + M_{N}(G_{2}) + 12\vert V_{2} \vert \vert E_{2} \vert^{2} M_{1}(G_{1}) + 8\vert E_{1} \vert \vert E_{2} \vert M_{1}(G_{2})\\
 && + 8\vert V_{2} \vert^{2} \vert E_{2} \vert (\vert V_{2} \vert + 1)M_{2}(G_{1}) + 8\vert V_{2} \vert \vert E_{1} \vert M_{2}(G_{2})\\
 && + 3\vert V_{2} \vert^{2}M_{1}(G_{1})M_{1}(G_{2}).
 \end{eqnarray*}
 Which is the desired result.

\begin{ex}The Wreath product of the path graphs $P_{n}$ and $P_{2}$ yield the Fence graph (figure. 6), whereas the Wreath product of the cycle $C_{n}$ and the path $P_{2}$ gives the closed Fence graph (figure. 6). Thus from the preposition \ref{pro4}, we compute the following results.
\begin{enumerate}
\item[(i)]
$M_{N}(P_{n}[P_{2}]) = 864n - 1694   ,n \geqslant 4$.
\item[(ii)]
$M_{N}(C_{n}[P_{2}]) = 816n + 2,     n \geqslant 3$.
\end{enumerate}
\end{ex}
\vspace{1cm}
\begin{figure}[h]
\begin{center}
\captionsetup{singlelinecheck = true}
\begin{tikzpicture}[scale=.4]
\tikzstyle{every node}=[draw, shape=circle, fill=black, scale=.4]
  \node (n1) at (3,1) {};
  \node (n2) at (3,4)  {};
  \node (n3) at (3,8)  {};
  \node (n4) at (3,11) {};
  \node (n5) at (3,14)  {};
  \node (n6) at (3,17)  {};
  \node (n7) at (9,1)  {};
  \node (n8) at (9,4)  {};
  \node (n9) at (9,8)  {};
  \node (n10) at (9,11)  {};
  \node (n11) at (9,14)  {};
  \node (n12) at (9,17)  {};
   \node (n13) at (6,5)  {};
  \node (n14) at (6,6)  {};
  \node (n15) at (6,7)  {};
  \foreach \from/\to in {n1/n2,n3/n4,n4/n5,n5/n6,n7/n8,n9/n10,n10/n11,n11/n12,n1/n7,n2/n8,n3/n9,n4/n10,n5/n11,n6/n12}
    \draw (\from) -- (\to);
\end{tikzpicture}
\hspace*{2.5cm}
\begin{tikzpicture}
\includegraphics[height=7cm,width=3.5cm]{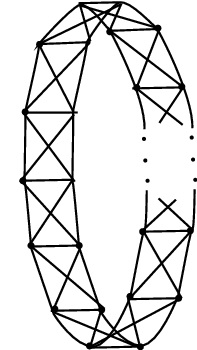}
\end{tikzpicture}
\caption*{~~~~~~~~~~~~~~~~~~~~~~~~~~~~~~~~~~~~~~$P_{n}[P_{2}]$ ~~~~~~~~~~~~~~~~~~~~~~~~~~~~~~~~~~~~~~~~~~~~~~~~~~~~~~~~~~~ $C_{n}[P_{2}]$}
\captionsetup{singlelinecheck = true, justification=justified}
\caption{Fence graph ({$P_{4}[P_{2}]$}) and closed Fence graph ($C_{n}[P_{2}]$).}
\end{center}
\end{figure}

\section{Conclusion}
In this article, we introduce the Neighbourhood Zagreb index($M_{N}$), investigate chemical applicability, compute some exact formulae for $M_{N}$ of some product graphs and apply the results on some chemical graphs. As a future work, we derive the results for some other graph operations, compute some bounds of this index and create some linear models with another indices having good correlation with different physiochemical properties of molecules. As the pharmacological activity of a compound depends on its physiochemical properties and the correlations of $M_{N}$ index with some of these properties are good, there is nothing to be surprised that $M_{N}$ index can be used in designing new drugs.

\section{Acknowledgements}
The first author is very obliged to the Department of Science and Technology (DST), Government of India for the Inspire Fellowship [IF170148].

\end{document}